\documentclass{birkmult}
\usepackage{amsmath}
\usepackage{amssymb,latexsym}
\usepackage{url}
\usepackage{graphicx}
\usepackage{longtable}
\usepackage{booktabs}
\usepackage{psfrag}

\newtheorem{deff}{Definition}
\newtheorem{thm}[deff]{Theorem}
\newtheorem{conj}[deff]{Conjecture}

\makeatletter

\title[Enumeration and Random Realization of Triangulated Surfaces]{Enumeration and Random Realization\\ of Triangulated Surfaces}

\author{Frank H.~Lutz}

\address{Technische Universit\"at Berlin\\
Fakult\"at II - Mathematik und Naturwissenschaften\\
Institut f\"ur Mathematik, Sekr.\ MA 3-2\\
Stra\ss e des 17.\ Juni 136\\
10623 Berlin\\
Germany}
\email{lutz@math.tu-berlin.de}

\date{}

\begin{document}

\begin{abstract}
We discuss different approaches for the enumeration of triangulated surfaces.
In particular, we enumerate all triangulated surfaces with 
$9$~and $10$ vertices. 
We also show how geometric realizations
of orientable surfaces with few vertices can be obtained 
by choosing coordinates randomly.
\end{abstract}

\maketitle

\section{Introduction}

The enumeration of triangulations of the $2$-dimensional sphere $S^2$
\index{triangulated surface!enumeration}
was started by Br\"uckner \cite{Brueckner1897} at the end of the 19th century.
It took Br\"uckner several years and nine thick manuscript books
to compose a list of triangulated $2$-spheres with up to $13$ vertices; cf.\ \cite{Brueckner1931}.
His enumeration was complete and correct up to $10$ vertices. On $11$, $12$, and
$13$ vertices his census comprised $1251$, $7616$, and $49451$ triangulations,
respectively. These numbers were off only slightly; they were later corrected by
Grace \cite{Grace1965} (for $11$ vertices), by Bowen and Fisk
\cite{BowenFisk1967}  (for $12$ vertices), and by Royle~\cite{Royle_url}
(for $13$ vertices). In fact, Royle made use of the
program \texttt{plantri} by Brinkmann and McKay~\cite{plantri}
to enumerate triangulations of the $2$-sphere with up to $23$ vertices
(see the manual of \texttt{plantri} and also Royle~\cite{Royle_url}).
Table~\ref{tbl:ten2d_spheres} and Table~\ref{tbl:ten2d_9_10} 
list the respective numbers.
Precise formulas for rooted triangulations of the $2$-sphere with $n$
vertices were determined by Tutte~\cite{Tutte1962}.

It follows from work of Steinitz \cite[\S 46]{SteinitzRademacher1934}
that every triangulated $2$-sphere can be reduced to the boundary
of the tetrahedron by a sequence of \emph{edge contractions}. In other words,
the boundary of the tetrahedron is the only \emph{irreducible
triangulation} of the $2$-sphere from which every $n$-vertex
\index{triangulated surface!irreducible}
triangulation can be obtained by a suitable sequence of \emph{vertex splits}.
The program \texttt{plantri} implements this procedure and allows for
a fast enumeration of triangulations of the $2$-sphere $S^2$.

Barnette and Edelson \cite{BarnetteEdelson1988} have shown
that every $2$-manifold has only finitely many irreducible triangulations.
The respective numbers were determined for the torus (21 examples) by
Gr\"unbaum and Lavrenchenko \cite{Lavrenchenko1990},
for the projective plane (2 examples) by Barnette \cite{Barnette1982b},
for the Klein bottle (29 examples) 
by Lawrencenko and Negami \cite{LawrencenkoNegami1997} 
and Sulanke~\cite{Sulanke2004pre},
and, recently, for the orientable surface of genus~2 (396784 examples)
and the non-orientable surfaces of genus~3 (9708 examples) and genus~4
(6297982 examples) by Sulanke \cite{Sulanke2005apre} \cite{Sulanke2005bpre}.
With his program \texttt{surftri} \cite{Sulanke2005cpre} (based on \texttt{plantri}),
Sulanke generated all triangulations with up to (at least) $14$ vertices 
for these surfaces by vertex-splitting; see \cite{Sulanke2005cpre} 
for respective counts of triangulations.

Here, we will survey further enumeration approaches.
In particular, we give a complete enumeration of all 
triangulated surfaces with up to 10 vertices.

{\begin{table}
\small\centering
\defaultaddspace=0.15em
\caption{Triangulated $2$-spheres with $11\leq n\leq 23$ vertices.}\label{tbl:ten2d_spheres}
\begin{tabular}{r@{\hspace{25mm}}r}\\
\toprule
 \addlinespace
 \addlinespace
 \addlinespace
  $n$ &   Types \\ \midrule
 \addlinespace
 \addlinespace
 \addlinespace
 \addlinespace
  11  &  1249 \\
 \addlinespace
  12  &  7595 \\
 \addlinespace
  13  &  49566 \\
 \addlinespace
  14  &  339722 \\
 \addlinespace
  15  &  2406841 \\
 \addlinespace
  16  &  17490241 \\
 \addlinespace
  17  &  129664753 \\
 \addlinespace
  18  &  977526957 \\
 \addlinespace
  19  &  7475907149 \\
 \addlinespace
  20  &  57896349553 \\
 \addlinespace
  21  &  453382272049 \\
 \addlinespace
  22  &  3585853662949 \\
 \addlinespace
  23  &  28615703421545 \\
 \addlinespace
 \addlinespace
 \addlinespace
\bottomrule
\end{tabular}
\end{table}
}

\bigskip

For an arbitrary orientable or non-orientable surface $M$ the Euler
characteristic $\chi(M)$ of $M$ equals, by Euler's equation,
the alternating sum of the number of vertices $n=f_0$, 
the number of edges $f_1$, and the number of triangles $f_2$, i.e.,
\begin{equation}
n-f_1+f_2=\chi (M).
\end{equation}
By double counting of incidences between edges and triangles of a triangulation, 
it follows that $2f_1=3f_2$. Thus, the number of vertices 
$n$ determines $f_1$ and $f_2$, that is, a triangulated surface $M$ 
of Euler characteristic $\chi (M)$ on $n$ vertices has \emph{$f$-vector}
\begin{equation}
f=(f_0,f_1,f_2)=(n,3n-3\chi(M),2n-2\chi(M)).
\end{equation}

\begin{figure}
\begin{center}
\psfrag{1}{\small 1}
\psfrag{2}{\small 2}
\psfrag{3}{\small 3}
\psfrag{4}{\small 4}
\psfrag{5}{\small 5}
\psfrag{6}{\small 6}
\psfrag{7}{\small 7}
\includegraphics[height=36mm]{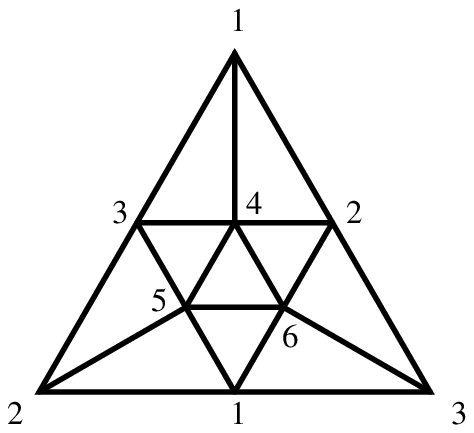}\hspace{15mm}\includegraphics[height=36mm]{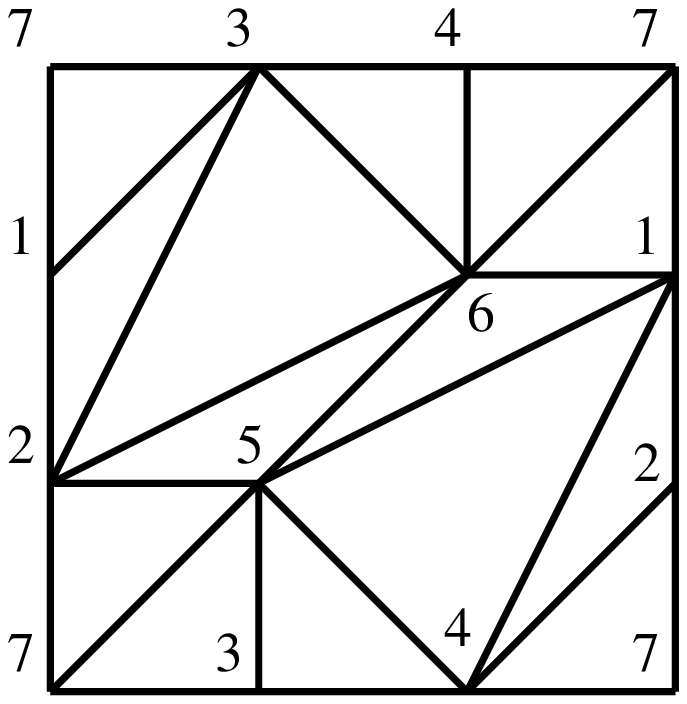}
\end{center}
\caption{The real projective plane ${\mathbb R}{\bf P}^2_6$ and M\"obius' torus.}
\label{fig:RP2_moebius}
\end{figure}

An \emph{orientable surface} $M(g,+)$ of genus $g$ has
homology $H_*(M(g,+))=({\mathbb Z},{\mathbb Z}^{2g},{\mathbb Z})$
and Euler characteristic $\chi(M(g,+))=2-2g$,
whereas a \emph{non-orient\-able surface} $M(g,-)$ of genus $g$ has
homology $H_*(M(g,-))=({\mathbb Z},{\mathbb Z}^{g-1}\oplus{\mathbb Z}_2,0)$
and Euler characteristic $\chi(M(g,-))=2-g$. 
The smallest possible $n$ for a triangulation of a $2$-manifold $M$
(with the exception of  the orientable
surface of genus~$2$, the Klein bottle, and the non-orient\-able
surface of genus~$3$, where an extra vertex has to be added, respectively)
\index{triangulated surface!vertex-minimal}
is determined by Heawood's bound~\cite{Heawood1890}
\begin{equation}
n\geq\Bigl\lceil\tfrac{1}{2}(7+\sqrt{49-24\chi (M)})\Bigl\rceil.
\end{equation}
Corresponding minimal and combinatorially unique triangulations 
of the real projective plane ${\mathbb R}{\bf P}^2$ with $6$ vertices (${\mathbb R}{\bf P}^2_6$)
and of the $2$-torus with $7$ vertices (M\"obius' torus \cite{Moebius1886})
were already known in the 19th century; see Figure~\ref{fig:RP2_moebius}.
However, it took until 1955 to complete the construction of series
of examples of minimal triangulations for
all non-orientable surfaces (Ringel~\cite{Ringel1955}) and 
until 1980 for all orientable surfaces (Jungerman and Ringel~\cite{JungermanRingel1980}).

A complete classification of triangulated surfaces with up to $8$ vertices
was obtained by Datta \cite{Datta1999}
and Datta and Nilakantan \cite{DattaNilakantan2002}.

\medskip
\smallskip

In Section~\ref{sec:enumeration}, we discuss algorithms
for the enumeration of triangulated surfaces,
and we give the numbers (up to combinatorial
equivalence) of triangulated surfaces
with $9$ and $10$ vertices.

For triangulated orientable surfaces of genus $g\geq 1$
it is, in general, a difficult problem to decide realizability.
Surprisingly, for triangulations with few vertices, $3$-dimensional
geometric realizations (with straight edges, flat triangles,
and without self-intersections) can be obtained by choosing
coordinates randomly; see Section~\ref{sec:random_realization}.

\section{Enumeration of Triangulated Surfaces}
\label{sec:enumeration}

At present, there are three essentially different enumeration schemes
\index{triangulated surfaces!enumeration}
available to algorithmically generate triangulated surfaces
with a fixed number $n$ of vertices:

\medskip

\begin{itemize}
\item[1.] \emph{Generation from irreducible triangulations}

          \vspace{1.5mm}

          \noindent
          An edge of a triangulated surface is \emph{contractible} 
          if the vertices of the edge can be identified without changing 
          the topological type of the triangulation.
          A triangulation of a surface is \emph{irreducible} if it has
          no contractible edge. 

          The triangulations of a given $2$-manifold $M$ 
          with $n$ vertices can therefore be obtained in two steps:
          First, determine all irreducible triangulations of $M$
          with $\leq n$ vertices, and, second, generate additional
          reducible triangulations of $M$ from the irreducible ones
          with $<n$ vertices by \emph{vertex-splitting}. 

          \vspace{1mm}
          
          \emph{Comments:} Although every $2$-manifold $M$ has only
          finitely many irreducible triangulations \cite{BarnetteEdelson1988},
          it is non-trivial to classify or enumerate these. In the
          case of $S^2$, however, the boundary of the tetrahedron
          is the only irreducible triangulation. Resulting reducible
          triangulations of $S^2$ with up to $23$ vertices were generated by
          Royle \cite{Royle_url} (with the program \texttt{plantri} 
          by Brinkmann and McKay~\cite{plantri}). For further results,
          in particular by Sulanke, see the overview in the introduction.

          \vspace{2mm}

\item[2.] \emph{Strongly connected enumeration}

          \vspace{1.5mm}

          \noindent
          The basic idea here is to start with a single triangle,
          which has three edges as its boundary. Then select the
          lexicographically smallest edge of the boundary. 
          This edge has to lie in a second
          \index{triangulated surface!enumeration!strongly connected}
          triangle (since in a triangulated surface every edge is contained
          in exactly two triangles). In order to complement our
          pivot edge to a triangle we can choose as a third vertex
          either a vertex of the current boundary or a new vertex
          that has not yet been used. For every such choice we 
          compose a new complex by adding the respective triangle
          to the previous complex. We continue with adding triangles 
          until, eventually, we obtain a closed surface.

          \vspace{1mm}

          \emph{Comments:} At every step of the procedure, the
          complexes that we produce are strongly connected.
          By the lexicographic choice of the pivot edges,
          the vertex-stars are closed in lexicographic order
          (which helps to sort out pseudomanifolds at an 
          early stage of the generation; see the discussion below).
          Nevertheless, the generation is not very systematic:
          Albeit we choose the pivot edges in lexicographic
          order the resulting triangulations do not have to
          be lexicographically minimal. 

          Strongly connected enumeration was used by
          Altshuler and Steinberg \cite{Altshuler1974},~\cite{AltshulerSteinberg1976}
          to determine all combinatorial $3$-manifolds with up to $9$
          vertices, by Bokowski \cite{AltshulerBokowskiSchuchert1996}
          (cf.\ also \cite{Bokowski2006pre}) to enumerate 
          all 59 neighborly triangulations with 12 vertices 
          of the orientable surface of genus $6$, and by
          Lutz and Sullivan to enumerate all 4787 triangulated $3$-manifolds
          of edge degree~$\leq 5$.

          \vspace{2mm}

\item[3.] \emph{Lexicographic enumeration}

          \vspace{1.5mm}

          \noindent
          It is often very useful to list a collection of
          objects, e.g., triangulated surfaces, in \emph{lexicographic order}:
          \index{triangulated surface!enumeration!lexicographic}
          Every listed triangulated surface with $n$ vertices
          is the lexicographically smallest set of triangles 
          combinatorially equivalent to this triangulation 
          and is lexicographically smaller than the next surface in
          the list.

          \vspace{2mm}


\item[3'.] \emph{Mixed lexicographic enumeration}

          \vspace{1.5mm}

          \noindent
          A variant of 3.; see below.
\end{itemize}

\bigskip

We present in the following an algorithm for the lexicographic
and the mixed lexicographic enumeration of triangulated surfaces.
(Similar enumeration schemes for vertex-transitive triangulations
have been described earlier by 
K\"uhnel and Lassmann \cite{KuehnelLassmann1985-di}
and K\"ohler and Lutz \cite{KoehlerLutz2005pre}.)

\bigskip

Let $\{1,2,\dots,n\}$ be the ground set of $n$ vertices.
Then a \emph{triangulation of a surface/a triangulated surface} with
$n$ vertices is

\smallskip

\begin{itemize}
\item a connected $2$-dimensional simplicial complex $M\subseteq 2^{\{1,\dots,n\}}$

      \vspace{1.5mm}

\item such that the link of every vertex of $M$ is a circle.
\end{itemize}

\smallskip

\noindent
In particular, 

\smallskip

\begin{itemize}
\item $M$ is \emph{pure}, that is, every maximal face of $M$ is $2$-dimensional,

      \vspace{1.5mm}

\item and every edge of $M$ is contained in exactly two triangles.
\end{itemize}

\smallskip

\noindent
As an example, we consider the case $n=5$.
On the ground set $\{1,2,\dots,5\}$ there are $\binom{5}{3}=10$ triangles,
$\binom{5}{2}=10$ edges, and $2^{10}$ different sets of triangles,
of which only few compose a triangulated surface with $5$ vertices.
Our aim will be to find those sets of triangles that indeed
form a triangulated surface. One way to proceed is by \emph{backtracking}:

\medskip

\begin{quote}
Start with some triangle and add further triangles
as long as no edge is contained in more than two
triangles. If this condition is violated, then backtrack.
A set of triangles is \emph{closed} if every of its edges
is contained in exactly two triangles. If the link of 
every vertex of a closed set of triangles is a circle,
then this set of triangles gives a triangulated surface: OUTPUT surface.
\end{quote}

\medskip

\noindent
We are interested in enumerating triangulated surfaces
\emph{up to combinatorial equivalence}, i.e., up to relabeling the vertices.
Thus we can, without loss of generality, assume that
the triangle $123$ should be present in the triangulation
and therefore can be chosen as the starting triangle.
More than that, we can assume that the triangulated surfaces 
that we are going to enumerate should come in lexicographic order.

In a lexicographically minimal triangulation, the collection $B_{{\rm deg}(1)}$ of triangles 
containing the vertex $1$ is of the form 
$$123,\, 124,\, 135,\, 146,\, 157,\, 168,\, \dots,\,  1\,({\rm deg}(1)-1)\,({\rm deg}(1)+1),\,
1\,{\rm deg}(1)\,({\rm deg}(1)+1),$$
where ${\rm deg}(1)$ is the \emph{degree} of the vertex $1$,  
i.e., the number of neighbors of the vertex~$1$. Obviously, 
$3\leq {\rm deg}(1)\leq n-1$.

\begin{figure}
  \begin{center}
\psfrag{1}{\small 1}
\psfrag{2}{\small 2}
\psfrag{3}{\small 3}
\psfrag{4}{\small 4}
\psfrag{5}{\small 5}
    \psfrag{d3}{\small $B_3$}
    \psfrag{d4}{\small $B_4$}
    \includegraphics[height=28mm]{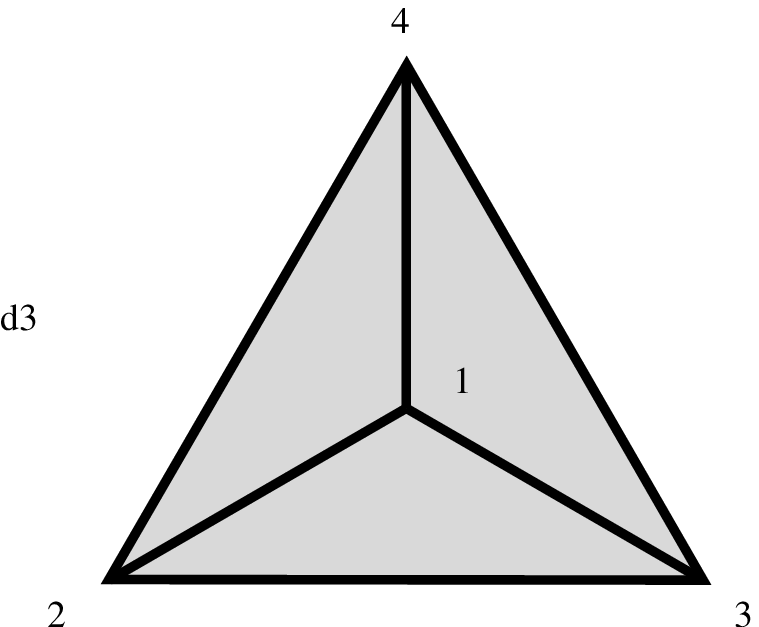}\hspace{25mm}\includegraphics[height=28mm]{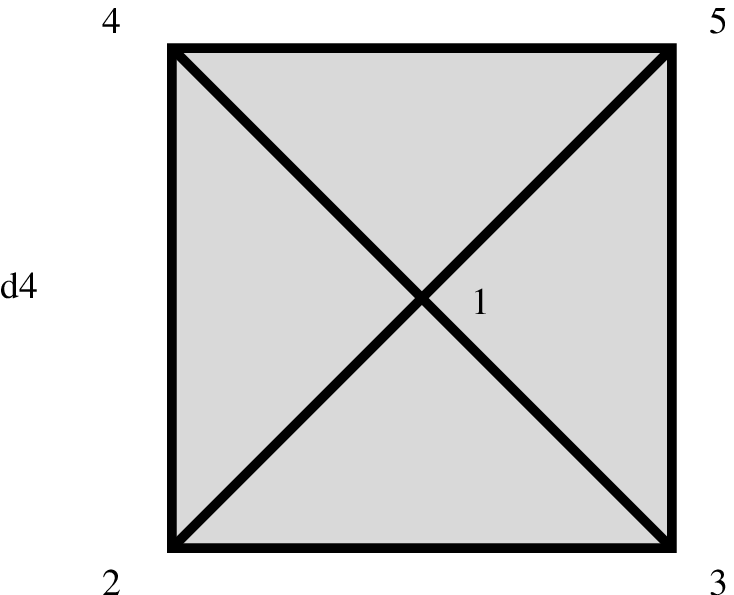}
  \end{center}
  \caption{Beginning segments}
  \label{fig:four}
\end{figure}

On $5$ vertices, the vertex $1$ has the possible beginning segments 
$B_3$ and $B_4$; see Figure~\ref{fig:four}.

In a lexicographically sorted list of (lexicographically smallest)
triangulated surfaces with $n$ vertices, those with beginning segment
$B_k$ are listed before those with beginning segment $B_{k+1}$, etc.
Thus, we start the backtracking with beginning segment $B_{3}$
and enumerate all corresponding triangulated surfaces,
then restart the backtracking with beginning segment $B_{4}$,
and so on. In triangulations with beginning segment $B_k$ 
all vertices have degree \emph{at least} $k$. Otherwise such a
triangulation has a vertex of degree \mbox{$j<k$}. However, by relabeling the vertices,
it can be achieved that the relabeled triangulation has beginning
segment $B_j$ (and thus would have appeared earlier in the
lexicographically sorted list). Contradiction.

\begin{figure}
  \addtolength{\arraycolsep}{.85mm} 
  $$\begin{array}{rl}
    \begin{array}{@{\extracolsep{10pt}}c@{\extracolsep{10pt}}c@{\extracolsep{10pt}}c@{\extracolsep{10pt}}c@{\extracolsep{10pt}}c@{\extracolsep{10pt}}c@{\extracolsep{10pt}}c@{\extracolsep{10pt}}c@{\extracolsep{10pt}}c@{\extracolsep{9pt}}c@{\extracolsep{14pt}}c} 12 & 13 & 14 & 15 & 23 & 24 & 25 & 34 & 35 & 45 \end{array} \\[.5mm]
    \begin{array}{c} 123\\124\\125\\134\\135\\145\\234\\235\\245\\345 \end{array} 
    \left(
      \begin{array}{cccccccccc}
        1 & 1 &   &   & 1 &   &   &   &   &   \\
        1 &   & 1 &   &   & 1 &   &   &   &   \\
        1 &   &   & 1 &   &   & 1 &   &   &   \\
          & 1 & 1 &   &   &   &   & 1 &   &   \\
          & 1 &   & 1 &   &   &   &   & 1 &   \\
          &   & 1 & 1 &   &   &   &   &   & 1 \\
          &   &   &   & 1 & 1 &   & 1 &   &   \\
          &   &   &   & 1 &   & 1 &   & 1 &   \\
          &   &   &   &   & 1 & 1 &   &   & 1 \\
          &   &   &   &   &   &   & 1 & 1 & 1 \\
      \end{array}
    \right) 
    & \begin{array}{c} \\ \\ \\ \end{array}
  \end{array}$$
  \addtolength{\arraycolsep}{-.85mm} 
\caption{The triangle-edge-incidence matrix for $n=5$ vertices.}
\label{fig:matrix}
\end{figure}

We store the triangles on the ground set as rows of a (sparse)
triangle-edge-incidence matrix; see Figure~\ref{fig:matrix}
for the triangle-edge-incidence matrix in the case of $n=5$ vertices.
The backtracking in terms of the triangle-edge-incidence matrix
then can be formulated as follows: Start with the zero row vector
and add to it all rows corresponding to the triangles of the
beginning segment $B_{3}$. The resulting vector has entries\\[3mm]
\begin{tabular}{@{\hspace{10pt}}ll}
    $0$ & (the corresponding edge does not appear in $B_{3}$), \\
    $1$ & (the corresponding edge is a boundary edge of $B_{3}$),\\ 
    $2$ & (the corresponding edge appears twice in the triangles of $B_{3}$). 
\end{tabular}\\[3mm]
We next add (the corresponding rows of) further triangles to (the sum
vector of) our beginning segment. As soon as a resulting entry is larger than
two, we backtrack, since in such a combination of triangles
the edge corresponding to the entry is contained in at least
three triangles, which is forbidden.
If a resulting vector has entries 0 and 2 only,
then the corresponding set of triangles is closed
and thus is a candidate for a triangulated surface.
In case of $n=5$ vertices, the backtracking (in short) 
is as follows:
    \addtolength{\arraycolsep}{-.72mm} 
  $$\begin{array}{l@{\hspace{5.25mm}}cccccccccccc@{\hspace{5.25mm}}l}
                   &   & 12 & 13 & 14 & 15 & 23 & 24 & 25 & 34 & 35 & 45\\[1.5mm]
  -:               & ( & 0 & 0 & 0 & 0 & 0 & 0 & 0 & 0 & 0 & 0 & ) & \\
 B_3:              & ( & 2 & 2 & 2 & 0 & 1 & 1 & 0 & 1 & 0 & 0 & ) & \\
 B_3+234:          & ( & 2 & 2 & 2 & 0 & 2 & 2 & 0 & 2 & 0 & 0 & ) & \mbox{\emph{Candidate!}} \\
 B_3:              & ( & 2 & 2 & 2 & 0 & 1 & 1 & 0 & 1 & 0 & 0 & ) & \\
 B_3+235:          & ( & 2 & 2 & 2 & 0 & 2 & 1 & 1 & 1 & 1 & 0 & ) & \\
 B_3+235+245:      & ( & 2 & 2 & 2 & 0 & 2 & 2 & 2 & 1 & 1 & 1 & ) & \\
 B_3+235+245+345:  & ( & 2 & 2 & 2 & 0 & 2 & 2 & 2 & 2 & 2 & 2 & ) & \mbox{\emph{Candidate!}} \\
 B_3+235+245:      & ( & 2 & 2 & 2 & 0 & 2 & 2 & 2 & 1 & 1 & 1 & ) & \\
 B_3+235:          & ( & 2 & 2 & 2 & 0 & 2 & 1 & 1 & 1 & 1 & 0 & ) & \\
 B_3:              & ( & 2 & 2 & 2 & 0 & 1 & 1 & 0 & 1 & 0 & 0 & ) & \\
  -:               & ( & 0 & 0 & 0 & 0 & 0 & 0 & 0 & 0 & 0 & 0 & ) & \\
 B_4:              & ( & 2 & 2 & 2 & 2 & 1 & 1 & 0 & 0 & 1 & 1 & ) &  \\
 B_4+234:          & ( & 2 & 2 & 2 & 2 & 2 & 2 & 0 & 1 & 1 & 1 & ) &  \\
 B_4+234+345:      & ( & 2 & 2 & 2 & 2 & 2 & 2 & 0 & 2 & 2 & 2 & ) & \mbox{\emph{Candidate!}} \\
 B_4+234:          & ( & 2 & 2 & 2 & 2 & 2 & 2 & 0 & 1 & 1 & 1 & ) & \\
 B_4:              & ( & 2 & 2 & 2 & 2 & 1 & 1 & 0 & 0 & 1 & 1 & ) & \\
 B_4+235:          & ( & 2 & 2 & 2 & 2 & 2 & 1 & 1 & 0 & 2 & 1 & ) & \\
 B_4+235+245:      & ( & 2 & 2 & 2 & 2 & 2 & 2 & 2 & 0 & 2 & 2 & ) & \mbox{\emph{Candidate!}} \\
 B_4+235:          & ( & 2 & 2 & 2 & 2 & 2 & 1 & 1 & 0 & 2 & 1 & ) & \\
 B_4:              & ( & 2 & 2 & 2 & 2 & 1 & 1 & 0 & 0 & 1 & 1 & ) & \\
  -:               & ( & 0 & 0 & 0 & 0 & 0 & 0 & 0 & 0 & 0 & 0 & ) & \\
  \end{array}$$
  \addtolength{\arraycolsep}{+.72mm} 
It remains to check whether every candidate is
indeed a triangulated surface. For this we have
to verify that the link of every vertex
is a circle, which is a purely combinatorial
condition that can easily be tested.
Moreover, we have to check for every candidate,
whether it is combinatorially isomorphic to
a triangulation that has appeared in the 
backtracking before. For example, the triangulations
$B_4+234+345$ and $B_4+235+245$ are isomorphic to $B_3+235+245+345$.
In fact, we should have stopped to add triangles to $B_4+234$: 
The link of vertex~$2$ in this case is a closed circle of length $3$,
which should not occur, since the beginning segment $B_4$ has degree $4$. 
Similarly, we should have backtracked at $B_4+235$.
Finally, we neglect triangulations with less than
$n$ vertices (such as $B_3+234$).
As result in the case $n=5$, we obtain that
there is a unique triangulation of the $2$-sphere
with $5$ vertices.

The basic property that every edge is contained in exactly two triangles is often
called \emph{the pseudomanifold property}: Every closed 
set of triangles forms a \emph{$2$-dimensional pseudomanifold}.

\begin{figure}
  \begin{center}
    \includegraphics[height=50mm]{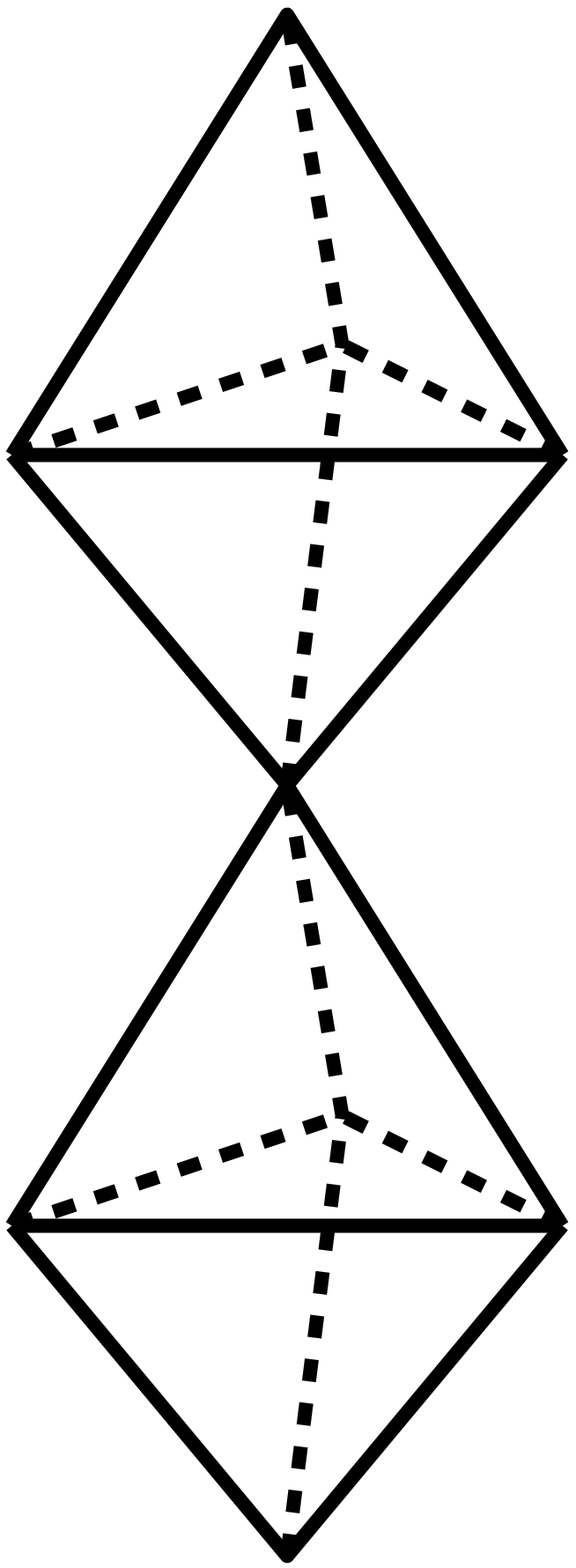}\hspace{45mm}\includegraphics[height=50mm]{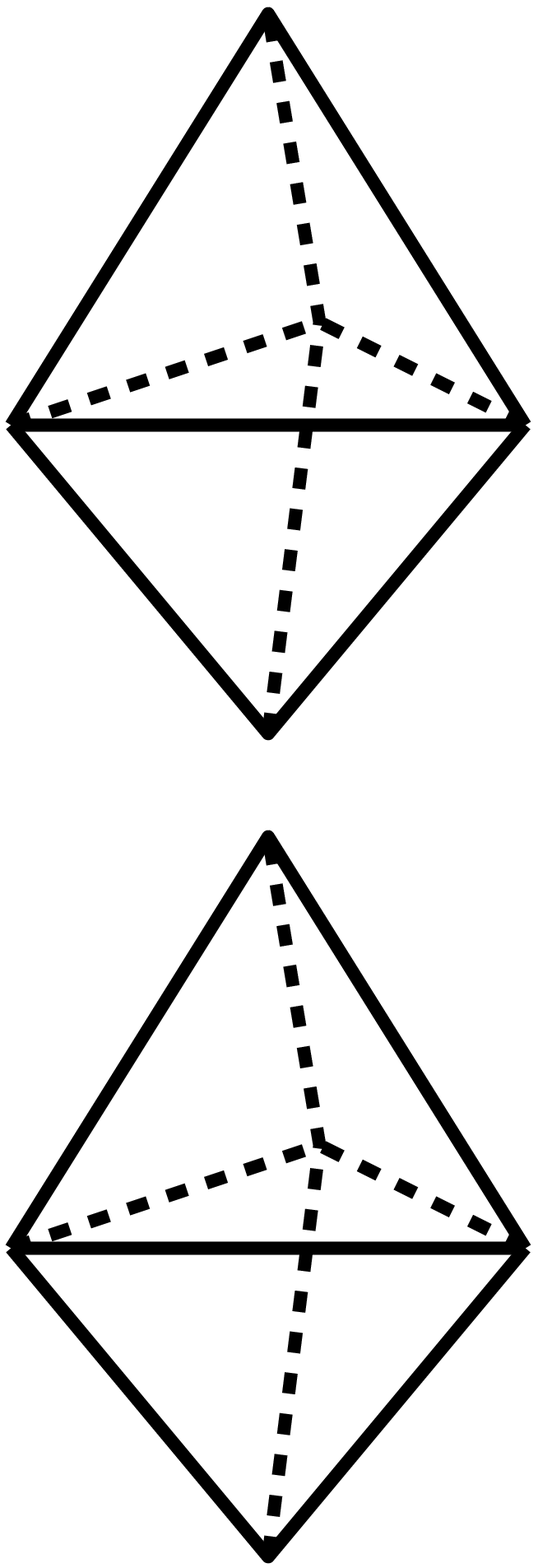}
  \end{center}
  \caption{Examples of pseudomanifolds}
  \label{fig:pseudo}
\end{figure}

It is rather easy to construct $2$-dimensional pseudomanifolds
that are not surfaces: The left pseudomanifold in Figure~\ref{fig:pseudo}
has the middle vertex as an isolated singularity with its vertex-link
consisting of two disjoint triangles. The right pseudomanifold in Figure~\ref{fig:pseudo}
has no singularities and thus is a surface, but it is not connected.

Since there are way more $2$-dimensional pseudomanifolds than there are 
(connected) surfaces, it is necessary to sort out those
pseudomanifolds that are not surfaces as early as possible in our backtracking.
In a $2$-dimensional pseudomanifold the link of every
vertex is a $1$-dimensional pseudomanifold (i.e., every vertex
lies in exactly two edges). Thus, every vertex-link
is a union of disjoint circles whereas the link of a vertex in a proper
surface is one single circle: 

\smallskip

\begin{itemize}
\item Backtrack if the link of a vertex of a current sum of triangles
  consists of a closed circle plus at least one extra edge.
\end{itemize}

\smallskip

Since at isolated singularities there are more edges in the link
than we want to have, we can try to avoid this as follows:

\smallskip

\begin{quote}
\emph{Mixed lexicographic enumeration:} Start the backtracking with
the beginning segment $B_{n-1}$ (instead of $B_3$), then proceed in
reversed order with $B_{n-2}$, \dots,
until $B_3$ is processed.
\end{quote}

\smallskip

\noindent
In the case that the beginning segment is $B_k$ 
we have the additional criterion:

\smallskip

\begin{itemize}
\item (Lexicographic enumeration:) Backtrack if the closed link of a
  vertex has less than $k$ edges.

\vspace{1mm}

\item (Mixed lexicographic enumeration:) Backtrack if the number of
  edges in the link of a vertex or if the degree of a vertex 
  is larger than $k$.
\end{itemize}

\smallskip

\begin{table}
\small\centering
\defaultaddspace=0.0em
\caption{Triangulated surfaces with up to $10$ vertices.}\label{tbl:ten2d_9_10}
\begin{tabular*}{\linewidth}{@{}l@{\extracolsep{12pt}}l@{\extracolsep{12pt}}r@{\extracolsep{\fill}}l@{\extracolsep{12pt}}l@{\extracolsep{12pt}}r@{\extracolsep{\fill}}l@{\extracolsep{12pt}}l@{\extracolsep{12pt}}r@{}}
\\\toprule\\[-3.5mm]
 \addlinespace
 \addlinespace
 \addlinespace
 \addlinespace
  $n$ &   Surface  &  Types  &    $n$ &   Surface  &  Types  &   $n$ &   Surface  &  Types   \\ \cmidrule{1-3}\cmidrule{4-6}\cmidrule{7-9}
 \addlinespace
 \addlinespace
 \addlinespace
 \addlinespace
  4  & $S^2$                  &   1 &   8  & $S^2$                  &  14 &   10  & $S^2$                  &   233 \\
 \addlinespace
     &                        &     &      & $T^2$                  &   7 &       & $T^2$                  &  2109 \\
 \addlinespace
  5  & $S^2$                  &   1 &      &                        &     &       & $M(2,+)$               &   865 \\
 \addlinespace
     &                        &     &      & ${\mathbb R}{\bf P}^2$ &  16 &       & $M(3,+)$               &    20 \\
 \addlinespace
  6  & $S^2$                  &   2 &      & $K^2$                  &   6 &       &                        &       \\
 \addlinespace
     &                        &     &      &                        &     &       & ${\mathbb R}{\bf P}^2$ &  1210 \\
 \addlinespace
     & ${\mathbb R}{\bf P}^2$ &   1 &   9  & $S^2$                  &  50 &       & $K^2$                  &  4462 \\
 \addlinespace
     &                        &     &      & $T^2$                  & 112 &       & $M(3,-)$               & 11784 \\
 \addlinespace
  7  & $S^2$                  &   5 &      &                        &     &       & $M(4,-)$               & 13657 \\
 \addlinespace
     & $T^2$                  &   1 &      & ${\mathbb R}{\bf P}^2$ & 134 &       & $M(5,-)$               &  7050 \\
 \addlinespace
     &                        &     &      & $K^2$                  & 187 &       & $M(6,-)$               &  1022 \\
 \addlinespace
     & ${\mathbb R}{\bf P}^2$ &   3 &      & $M(3,-)$               & 133 &       & $M(7,-)$               &    14 \\
 \addlinespace
     &                        &     &      & $M(4,-)$               &  37 \\
 \addlinespace
     &                        &     &      & $M(5,-)$               &   2 \\
 \addlinespace
 \addlinespace
 \addlinespace
 \addlinespace
\bottomrule
\end{tabular*}
\end{table}

In the mixed lexicographic approach, the sub-collections of triangulations
that contain a vertex of the same maximal degree are sorted lexicographically. 
In particular, every triangulation in such a sub-collection begins with the same
beginning segment (which is of type $B_k$ for some $k$).
However, the complete list of triangulations is not produced
in lexicographic order anymore.

By symmetry, we can, in addition, exclude the following cases
(for the lexicographic as well as for the mixed lexicographic approach):

\smallskip

\begin{itemize}
\item Do not use triangles of the form $23j$ with odd $5\leq j\leq k$.
  (Since a resulting surface would, by relabeling, be isomorphic
  to a triangulation with beginning segment $B_k$ plus triangle
  $23(j-1)$, which is lexicographically smaller).
  If the triangle $23i$ with even $6\leq i\leq k$ is used, then do not
  use the triangles $24j$ with odd $3\leq j\leq i-3$.
\end{itemize}

\smallskip

Finally, we test for every resulting (connected) surface whether
it has, up to combinatorial equivalence, appeared previously
in the enumeration. For this, we first compute as combinatorial invariants
the $f$-vector, the sequence of vertex degrees,
and the \emph{Altshuler-Steinberg determinant}~\cite{AltshulerSteinberg1973} 
(i.e., the determinant ${\rm det}(AA^T)$ of the vertex-triangle incidence matrix~$A$)
of an example. If these invariants are equal for two resulting surfaces,
then we take one triangle of the first complex and test for all
possible ways it can be mapped to the triangles of the second complex
whether this map can be extended to a simplicial isomorphism of the
two complexes. (Alternatively, one can use McKay's fast graph isomorphism testing program
\texttt{nauty} \cite{nauty} to determine whether the
vertex-facet incidence graphs of the two complexes are isomorphic or not.)

\begin{thm}
There are exactly $655$ triangulated surfaces with $9$ vertices
and $42426$ triangulated surfaces with $10$ vertices.
\end{thm}

The respective triangulations (in mixed lexicographic order) 
can be found online at \cite{Lutz_PAGE}. 
Table~\ref{tbl:ten2d_9_10} gives the detailed numbers
of orientable and non-orientable surfaces that appear
with up to $10$ vertices. The corresponding combinatorial
symmetry groups $G$ of the examples are listed in
Table~\ref{tbl:2man}.


\pagebreak

{\small
\defaultaddspace=-.005em

\setlength{\LTleft}{0pt}
\setlength{\LTright}{0pt}
\begin{longtable}{@{}l@{\extracolsep{8pt}}l@{\extracolsep{8pt}}r@{\extracolsep{12pt}}l@{\extracolsep{8pt}}r@{\extracolsep{\fill}}l@{\extracolsep{8pt}}l@{\extracolsep{8pt}}r@{\extracolsep{12pt}}l@{\extracolsep{8pt}}r@{}}
\caption{\protect\parbox[t]{9cm}{Symmetry groups of triangulated surfaces with up to $10$ vertices.}}\label{tbl:2man}
\\\toprule\\[-3mm]
 \addlinespace
 \addlinespace
 \addlinespace
 \addlinespace
 $n$ &   Surface              & $|G|$ &   $G$                         & Types & $n$ &   Surface  &  $|G|$  &   $G$ &  Types   \\ \cmidrule{1-5}\cmidrule{6-10}
\endfirsthead
\caption{\protect\parbox[t]{9cm}{Symmetry groups of triangulated surfaces (continued).}}
\\\toprule\\[-3mm]
 \addlinespace
 \addlinespace
 \addlinespace
 \addlinespace
 $n$ &   Surface              & $|G|$ &   $G$                         & Types & $n$ &   Surface  &  $|G|$  &   $G$ &  Types   \\ \cmidrule{1-5}\cmidrule{6-10}
\endhead
\bottomrule
\endfoot
 \addlinespace
  4  & $S^2$                  &  24 & $T^*=S_4$,                          &   &   9 & $S^2$                  &   1 & trivial                             & 16 \\
 \addlinespace
     &                        &     & $3$-transitive                      & 1 &     &                        &   2 & ${\mathbb Z}_2$                     & 25 \\
 \addlinespace
     &                        &     &                                     &   &     &                        &   4 & ${\mathbb Z}_2\times{\mathbb Z}_2$ & 5 \\
 \addlinespace
  5  & $S^2$                  &  12 & $S_3\times {\mathbb Z}_2$           & 1 &     &                        &   6 & $S_3$                               & 1 \\
 \addlinespace
     &                        &     &                                     &   &     &                        &  12 & $S_3\times {\mathbb Z}_2$           & 2 \\
 \addlinespace
  6  & $S^2$                  &   4 & ${\mathbb Z}_2\times {\mathbb Z}_2$ & 1 &     &                        &  28 & $D_7\times {\mathbb Z}_2$           & 1 \\
 \addlinespace
     &                        &  48 & $O^*={\mathbb Z}_2\wr S_3$,         &   &     & $T^2$                  &   1 & trivial                             & 52 \\
 \addlinespace
     &                        &     & vertex-trans.                       & 1 &     &                        &   2 & ${\mathbb Z}_2$                     & 46 \\
 \addlinespace
     & ${\mathbb R}{\bf P}^2$ &  60 & $A_5$,                              &   &     &                        &   3 & ${\mathbb Z}_3$                     & 1 \\
 \addlinespace
     &                        &     & $2$-transitive                      & 1 &     &                        &   4 & ${\mathbb Z}_2\times{\mathbb Z}_2$ & 7 \\
 \addlinespace
     &                        &     &                                     &   &     &                        &   6 & ${\mathbb Z}_6$                     & 1 \\
 \addlinespace
  7  & $S^2$                  &   2 & ${\mathbb Z}_2$                     & 1 &     &                        &     & $S_3$                               & 1 \\
 \addlinespace
     &                        &   4 & ${\mathbb Z}_2\times {\mathbb Z}_2$ & 1 &     &                        &  12 & $S_3\times {\mathbb Z}_2$           & 2 \\
 \addlinespace
     &                        &   6 & $S_3$                               & 2 &     &                        &  18 & $D_9$,                              &   \\
 \addlinespace
     &                        &  20 & $D_{10}$                            & 1 &     &                        &     & vertex-trans.                       & 1 \\
 \addlinespace
     & $T^2$                  &  42 & $AGL(1,7)$,                         &   &     &                        & 108 & ${\mathbb Z}_3^{\,2}\!:\!D_6$,      &   \\
 \addlinespace
     &                        &     & $2$-transitive                      & 1 &     &                        &     & vertex-trans.                       & 1 \\
 \addlinespace
     & ${\mathbb R}{\bf P}^2$ &   4 & ${\mathbb Z}_2\times {\mathbb Z}_2$ & 1 &     & ${\mathbb R}{\bf P}^2$ &   1 & trivial                             & 63 \\
 \addlinespace
     &                        &   6 & $S_3$                               & 1 &     &                        &   2 & ${\mathbb Z}_2$                     & 52 \\
 \addlinespace
     &                        &  24 & $T^*=S_4$                           & 1 &     &                        &   3 & ${\mathbb Z}_3$                     & 1 \\
 \addlinespace
     &                        &     &                                     &   &     &                        &   4 & ${\mathbb Z}_4$                     & 1 \\
 \addlinespace
  8  & $S^2$                  &   1 & trivial                             & 2 &     &                        &     & ${\mathbb Z}_2\times {\mathbb Z}_2$ & 11 \\
 \addlinespace
     &                        &   2 & ${\mathbb Z}_2$                     & 5 &     &                        &   6 & $S_3$                               & 5 \\
 \addlinespace
     &                        &   4 & ${\mathbb Z}_2\times {\mathbb Z}_2$ & 3 &     &                        &  24 & $T^*=S_4$                           & 1 \\
 \addlinespace
     &                        &   8 & $D_4$                               & 1 &     & $K^2$                  &   1 & trivial                             & 131 \\
 \addlinespace
     &                        &  12 & $S_3\times {\mathbb Z}_2$           & 1 &     &                        &   2 & ${\mathbb Z}_2$                     & 46 \\
 \addlinespace
     &                        &  24 & $T^*=S_4$                           & 1 &     &                        &   4 & ${\mathbb Z}_2\times {\mathbb Z}_2$ & 6 \\
 \addlinespace
     &                        &     & $D_6\times {\mathbb Z}_2$           & 1 &     &                        &   6 & $S_3$                               & 2 \\
 \addlinespace
     & $T^2$                  &   2 & ${\mathbb Z}_2$                     & 2 &     &                        &  12 & $S_3\times {\mathbb Z}_2$           & 2 \\
 \addlinespace
     &                        &   3 & ${\mathbb Z}_3$                     & 1 &     & $M(3,-)$               &   1 & trivial                             & 106 \\
 \addlinespace
     &                        &   4 & ${\mathbb Z}_2\times {\mathbb Z}_2$ & 2 &     &                        &   2 & ${\mathbb Z}_2$                     & 24 \\
 \addlinespace
     &                        &   6 & $S_3$                               & 1 &     &                        &   3 & ${\mathbb Z}_3$                     & 1 \\
 \addlinespace
     &                        &  32 & $[32,43]$,                          &   &     &                        &   6 & ${\mathbb Z}_6$                     & 1 \\
 \addlinespace
     &                        &     & vertex-trans.                       & 1 &     &                        &     & $S_3$                               & 1 \\
 \addlinespace
     & ${\mathbb R}{\bf P}^2$ &   1 & trivial                             & 2 &     & $M(4,-)$               &   1 & trivial                             & 19 \\
 \addlinespace
     &                        &   2 & ${\mathbb Z}_2$                     & 8 &     &                        &   2 & ${\mathbb Z}_2$                     & 11 \\
 \addlinespace
     &                        &   4 & ${\mathbb Z}_2\times {\mathbb Z}_2$ & 5 &     &                        &   3 & ${\mathbb Z}_3$                     & 1 \\
 \addlinespace
     &                        &  14 & $D_7$                               & 1 &     &                        &   4 & ${\mathbb Z}_2\times {\mathbb Z}_2$ & 5 \\
 \addlinespace
     & $K^2$                  &   1 & trivial                             & 1 &     &                        &   6 & $S_3$                               & 1 \\
 \addlinespace
     &                        &   2 & ${\mathbb Z}_2$                     & 4 &     & $M(5,-)$               &   6 & ${\mathbb Z}_6$                     & 1 \\
 \addlinespace
     &                        &   8 & $D_4$                               & 1 &     &                        &  18 & $S_3\times {\mathbb Z}_3$,          &   \\
 \addlinespace
     &                        &     &                                     &      &   &                        &     & vertex-trans.                      & 1 \\
 \addlinespace
 \addlinespace
 \addlinespace
 \addlinespace
\newpage
 \addlinespace
 10  & $S^2$                  &   1 & trivial                             &  137 &   & $K^2$                  &   1 & trivial                             & 4057 \\
 \addlinespace
     &                        &   2 & ${\mathbb Z}_2$                     &   69 &   &                        &   2 & ${\mathbb Z}_2$                     & 367 \\
 \addlinespace
     &                        &   3 & ${\mathbb Z}_3$                     &    1 &   &                        &   4 & ${\mathbb Z}_4$                     & 1 \\
 \addlinespace
     &                        &   4 & ${\mathbb Z}_2\times {\mathbb Z}_2$ &   13 &   &                        &     & ${\mathbb Z}_2\times {\mathbb Z}_2$ & 30 \\
 \addlinespace
     &                        &   6 & $S_3$                               &    6 &   &                        &   8 & ${\mathbb Z}_2^{\,3}$               & 2 \\
 \addlinespace
     &                        &   8 & ${\mathbb Z}_2^{\,3}$               &    2 &   &                        &     & $D_4$                               & 3 \\
 \addlinespace
     &                        &     & $D_4$                               &    2 &   &                        &  10 & $D_5$                               & 1 \\
 \addlinespace
     &                        &  16 & $D_8$                               &    1 &   &                        &  20 & $D_{10}$,                           &   \\
 \addlinespace
     &                        &  24 & $T^*=S_4$                           &    1 &   &                        &     & vertex-trans.                       & 1 \\
 \addlinespace
     &                        &  32 & $D_8\times {\mathbb Z}_2$           &    1 &   & $M(3,-)$               &   1 & trivial                             & 11308 \\
 \addlinespace
     & $T^2$                  &   1 & trivial                             & 1763 &   &                        &   2 & ${\mathbb Z}_2$                     & 448 \\
 \addlinespace
     &                        &   2 & ${\mathbb Z}_2$                     &  292 &   &                        &   3 & ${\mathbb Z}_3$                     & 12 \\
 \addlinespace
     &                        &   3 & ${\mathbb Z}_3$                     &   10 &   &                        &   4 & ${\mathbb Z}_2\times {\mathbb Z}_2$ & 11 \\
 \addlinespace
     &                        &   4 & ${\mathbb Z}_4$                     &    2 &   &                        &   6 & $S_3$                               & 5 \\
 \addlinespace
     &                        &     & ${\mathbb Z}_2\times {\mathbb Z}_2$ &   33 &   & $M(4,-)$               &   1 & trivial                             & 13037 \\
 \addlinespace
     &                        &   6 & $S_3$                               &    3 &   &                        &   2 & ${\mathbb Z}_2$                     & 556 \\
 \addlinespace
     &                        &   8 & ${\mathbb Z}_2^{\,3}$               &    2 &   &                        &   3 & ${\mathbb Z}_3$                     & 6 \\
 \addlinespace
     &                        &     & $D_4$                               &    1 &   &                        &   4 & ${\mathbb Z}_4$                     & 1 \\
 \addlinespace
     &                        &  12 & $S_3\times {\mathbb Z}_2$           &    1 &   &                        &     & ${\mathbb Z}_2\times {\mathbb Z}_2$ & 45 \\
 \addlinespace
     &                        &  20 & $D_{10}$,                           &      &   &                        &   6 & ${\mathbb Z}_6$                     & 1 \\
 \addlinespace
     &                        &     & vertex-trans.                       &    1 &   &                        &     & $S_3$                               & 5 \\
 \addlinespace
     &                        &     & $AGL(1,5)$                          &    1 &   &                        &   8 & ${\mathbb Z}_2^{\,3}$               & 1 \\
 \addlinespace
     & $M(2,+)$               &   1 & trivial                             &  789 &   &                        &     & $D_4$                               & 4 \\
 \addlinespace
     &                        &   2 & ${\mathbb Z}_2$                     &   61 &   &                        &  48 & $O^*={\mathbb Z}_2\wr S_3$            & 1 \\
 \addlinespace
     &                        &   3 & ${\mathbb Z}_3$                     &    7 &   & $M(5,-)$               &   1 & trivial                             & 6792 \\
 \addlinespace
     &                        &   4 & ${\mathbb Z}_4$                     &    2 &   &                        &   2 & ${\mathbb Z}_2$                     & 214 \\
 \addlinespace
     &                        &     & ${\mathbb Z}_2\times {\mathbb Z}_2$ &    4 &   &                        &   3 & ${\mathbb Z}_3$                     & 22 \\
 \addlinespace
     &                        &   8 & ${\mathbb Z}_8$                     &    1 &   &                        &   4 & ${\mathbb Z}_2\times {\mathbb Z}_2$ & 16 \\
 \addlinespace
     &                        &  16 & $\langle\, 2,2\,|\,2\,\rangle$      &    1 &   &                        &   6 & ${\mathbb Z}_6$                     & 2 \\
 \addlinespace
     & $M(3,+)$               &   1 & trivial                             &    8 &   &                        &     & $S_3$                               & 4 \\
 \addlinespace
     &                        &   2 & ${\mathbb Z}_2$                     &    3 &   & $M(6,-)$               &   1 & trivial                             & 926 \\
 \addlinespace
     &                        &   3 & ${\mathbb Z}_3$                     &    4 &   &                        &   2 & ${\mathbb Z}_2$                     & 71 \\
 \addlinespace
     &                        &   4 & ${\mathbb Z}_4$                     &    2 &   &                        &   3 & ${\mathbb Z}_3$                     & 18 \\
 \addlinespace
     &                        &     & ${\mathbb Z}_2\times {\mathbb Z}_2$ &    1 &   &                        &   4 & ${\mathbb Z}_2\times {\mathbb Z}_2$ & 6 \\
 \addlinespace
     &                        &  12 & $A_4$                               &    1 &   &                        &  12 & $A_4$                               & 1 \\
 \addlinespace
     &                        &  21 & ${\mathbb Z}_7\times {\mathbb Z}_3$ &    1 &   & $M(7,-)$               &   1 & trivial                             & 4 \\
 \addlinespace
     & ${\mathbb R}{\bf P}^2$ &   1 & trivial                             &  923 &   &                        &   2 & ${\mathbb Z}_2$                     & 4 \\
 \addlinespace
     &                        &   2 & ${\mathbb Z}_2$                     &  242 &   &                        &   3 & ${\mathbb Z}_3$                     & 1 \\
 \addlinespace
     &                        &   3 & ${\mathbb Z}_3$                     &    2 &   &                        &   5 & ${\mathbb Z}_5$                     & 1 \\
 \addlinespace
     &                        &   4 & ${\mathbb Z}_2\times {\mathbb Z}_2$ &   29 &   &                        &   6 & $S_3$                               & 1 \\
 \addlinespace
     &                        &   6 & $S_3$                               &   10 &   &                        &   9 & ${\mathbb Z}_9$                     & 1 \\
 \addlinespace
     &                        &  12 & $S_3\times {\mathbb Z}_2$           &    2 &   &                        &  12 & $A_4$                               & 1 \\
 \addlinespace
     &                        &     & $A_4$                               &    1 &   &                        &  60 & $A_5$,                              &   \\
 \addlinespace
     &                        &  18 & $D_9$                               &    1 &   &                        &     & vertex-trans.                       & 1 \\
 \addlinespace
 \addlinespace
 \addlinespace
 \addlinespace
\end{longtable}

}


\section{Random Realization}
\label{sec:random_realization}

It was asked by Gr\"unbaum~\cite[Ch.~13.2]{Gruenbaum1967},
whether every triangulated orientable surface can be embedded 
geometrically in\, ${\mathbb R}^3$, i.e., whether it can be realized 
\index{triangulated surface!realization}
with straight edges, flat triangles, and without self intersections?
By Steinitz' theorem (cf.~\cite{Steinitz1922},
\cite{SteinitzRademacher1934}, \cite{Ziegler1995}), every triangulated $2$-sphere 
is realizable as the boundary complex of a convex $3$-dimensional polytope. 
For the $2$-torus of genus $1$ the realizability problem is still open.
\begin{conj} {\rm (Duke~\cite{Duke1970})}
Every triangulated torus is realizable in ${\mathbb R}^3$.
\end{conj}

\begin{figure}
\begin{center}
\psfrag{1}{\small 1}
\psfrag{2}{\small 2}
\psfrag{3}{\small 3}
\psfrag{4}{\small 4}
\psfrag{5}{\small 5}
\psfrag{6}{\small 6}
\psfrag{7}{\small 7}
\psfrag{Coordinates:}{\small Coordinates:}
\psfrag{1:   (3,-3,0)}{\small 1:   $(3,-3,0)$}
\psfrag{2:   (-3,3,0)}{\small 2:   $(-3,3,0)$}
\psfrag{3:   (-3,-3,1)}{\small 3:   $(-3,-3,1)$}
\psfrag{4:   (3,3,1)}{\small 4:   $(3,3,1)$}
\psfrag{5:   (-1,-2,3)}{\small 5:   $(-1,-2,3)$}
\psfrag{6:   (1,2,3)}{\small 6:   $(1,2,3)$}
\psfrag{7:   (0,0,15)}{\small 7:   $(0,0,15)$}
\psfrag{Triangles:}{\small Triangles:}
\psfrag{123      145      156      345      167      467      247}{\small 123\,\, 145\,\, 156\,\, 345\,\, 167\,\, 467\,\, 247}
\psfrag{124      236      256      346      257      357      137}{\small 124\,\, 236\,\, 256\,\, 346\,\, 257\,\, 357\,\, 137}
\includegraphics[width=.4125\linewidth]{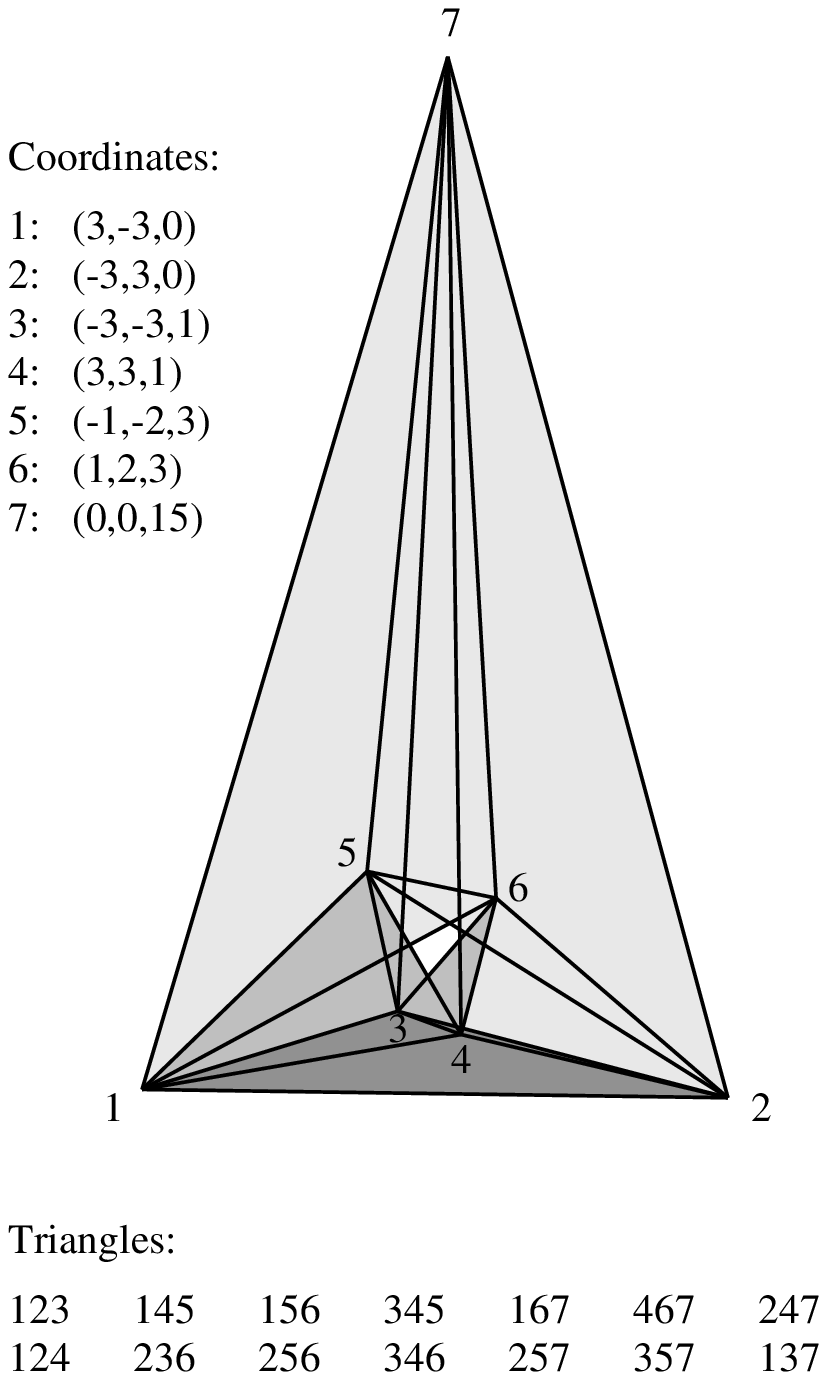}
\end{center}
\caption{{Cs\'asz\'ar}'s torus.}
\label{fig:7torus_small_ten2d}
\end{figure}

A first explicit geometric realization (see Figure~\ref{fig:7torus_small_ten2d}) 
of M\"obius' minimal $7$-vertex triangulation \cite{Moebius1886} of the $2$-torus 
was given by Cs\'asz\'ar~\cite{Csaszar1949} 
(cf.\ also \cite{Gardner1975b} and \cite{Lutz2002b}). Bokowski and Eggert
\cite{BokowskiEggert1991} showed that there are altogether $72$
different ``types'' of realizations of the M\"obius torus,
and Fendrich \cite{Fendrich2003} verified that 
triangulated tori with up to $11$ vertices are realizable.

Brehm and Bokowski \cite{BokowskiBrehm1987},
\cite{BokowskiBrehm1989}, \cite{Brehm1981}, \cite{Brehm1987b}
constructed geometric realizations for several triangulated 
orientable surfaces of genus\, $g=2,3,4$\,
with minimal numbers of vertices $n=10,10,11$, respectively.

Neighborly triangulations (i.e., triangulations that have as its
\index{triangulated surface!neighborly}
$1$-skeleton the complete graph $K_n$) of higher genus were considered
as candidates for counter-examples  to the Gr\"un\-baum realization
problem for a while (cf.~\cite{Csaszar1949} and \cite[p.~137]{BokowskiSturmfels1989}). 
Neighborly orientable surfaces have genus\, $g=(n-3)(n-4)/12$\, and
therefore\, $n\equiv 0,3,4,7\,{\rm mod}\,12$\, vertices, with\, $g=6$\, and\, $n=12$\, as the first
case beyond the tetrahedron and the $7$-vertex torus.

\begin{thm} {\rm (Bokowski and Guedes de Oliveira~\cite{BokowskiGuedes_de_Oliveira2000})}
The triangulated orientable surface $N^{12}_{54}$ of genus $6$ with $12$
vertices of Alshuler's list~\cite{AltshulerBokowskiSchuchert1996}
is not geometrically embeddable in ${\mathbb R}^3$.
\end{thm}

In fact, Bokowski and Guedes de Oliveira showed that there is no
oriented matroid compatible with the triangulation $N^{12}_{54}$,
from which the non-realizability follows.
Recently, Schewe reimplemented the approach of Bokowski and Guedes de Oliveira:

\begin{thm} {\rm (Schewe \cite{Schewe2006pre})}
Every orientable surface of genus $g\geq 5$ has a triangulation 
which is not geometrically embeddable in ${\mathbb R}^3$.
\end{thm}

From an algorithmic point of view the realizability problem
for triangulated surfaces is decidable 
(cf.~\cite[p.~50]{BokowskiSturmfels1989}),
but there is \emph{no} algorithm known that would
solve the realization problem for instances with, say, $10$ vertices
in reasonable time.

Surprisingly, the following simple heuristic can be used
to realize tori and surfaces of genus $2$ with up to $10$ vertices.

\bigskip

\noindent
\textbf{Random Realization}\, \emph{For a given orientable surface
with $n$ vertices pick $n$ integer points in a cube 
of size $k^3$ (for some fixed $k\in {\mathbb N}$)
uniformly at random. Test whether this (labeled) set
of $n$ points in ${\mathbb R}^3$ yields a geometric
realization of the surface. If not, try again.}
\index{triangulated surface!realization!random}

\bigskip

\noindent
Suppose, we are given $n$ (integer) points $\vec{x}_1,\dots,\vec{x}_{n}$ 
(in general position) in ${\mathbb R}^3$
together with a triangulation of an orientable surface.
It then is an elementary linear algebra exercise to check 
whether these $n$ points provide a geometric realization of the surface: 
For every pair of a triangle $i_1i_2i_3$ together with a combinatorially 
disjoint edge $i_4i_5$ of the triangulation we have to test whether 
the geometric triangle $\vec{x}_{i_1}\vec{x}_{i_2}\vec{x}_{i_3}$
and the edge $\vec{x}_{i_4}\vec{x}_{i_5}$ have empty intersection.

\bigskip

Since the surfaces were enumerated in (mixed) lexicographic order with
different triangulations in the list sometimes differing only slightly,
it is promising to try:

\bigskip

\noindent
\textbf{Recycling of Coordinates}\, \emph{As soon as 
a realization has been found for a triangulated orientable surface
with $n$ vertices, test whether the corresponding
set of coordinates yields realizations for other
triangulations with $n$ vertices as well. Moreover,
perturb the coordinates slightly and try again.}

\bigskip

\begin{figure}
\begin{center}
\psfrag{1}{\small 1}
\psfrag{2}{\small 2}
\psfrag{3}{\small 3}
\psfrag{4}{\small 4}
\psfrag{5}{\small 5}
\psfrag{6}{\small 6}
\psfrag{7}{\small 7}
\psfrag{8}{\small 8}
\psfrag{9}{\small 9}
\psfrag{10}{\small 10}
\includegraphics[width=.5\linewidth]{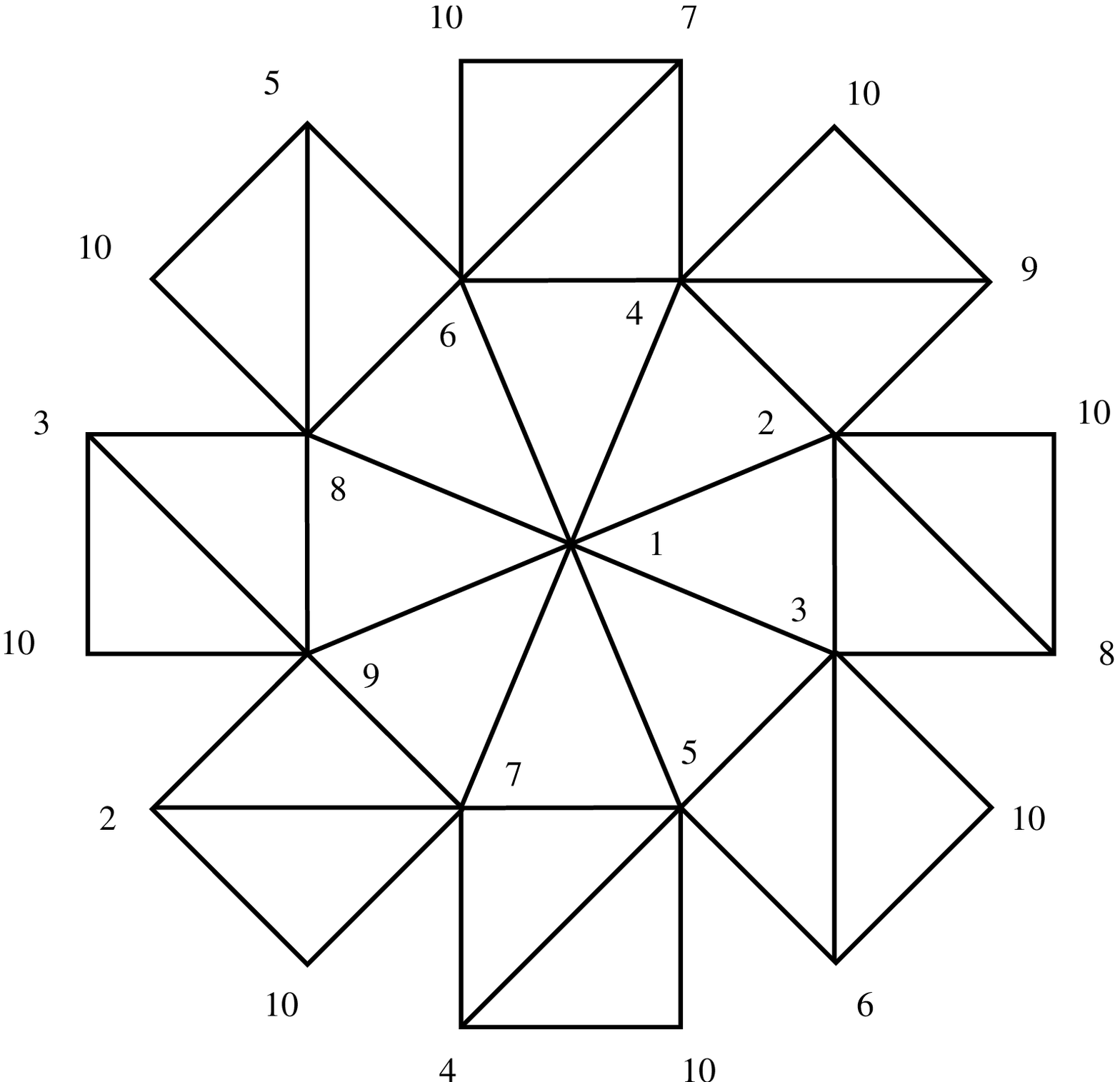}
\end{center}
\caption{The double torus\, \texttt{manifold\_2\_10\_41348}.}
\label{fig:manifold_41348}
\end{figure}

Random realization and recycling of coordinates was used
to obtain realizations for $864$ of the $865$ triangulations 
of the orientable surface of genus $2$ with $10$ vertices, 
leaving one case open. This last example (\texttt{manifold\_2\_10\_41348}
in the catalog \cite{Lutz_PAGE}) has the highest symmetry group
of order 16 among the $865$ examples; see Figure~\ref{fig:manifold_41348}.
It was realized geometrically by J\"urgen Bokowski 
by construction of an explicit rubber-band model.

\begin{thm} {\rm (Bokowski and Lutz; cf.\ \cite{Bokowski2006pre})}
All $865$ vertex-minimal $10$-vertex triangulations of the orientable surface 
of genus $2$ can be realized geometrically in~${\mathbb R}^3$.
\end{thm}
\begin{conj}
Every triangulation of the orientable surface of genus $2$ 
is realizable in ${\mathbb R}^3$.
\end{conj}
A priori, nothing is known on the expected number of tries
that are necessary to obtain a geometric realization for a given
triangulated orientable surface by random realization:
For triangulations of a surface of genus $g\geq 1$ 
any proof showing the finiteness of the expected number of tries
would immediately imply the realizability of the given triangulation.

\bigskip

\noindent
\emph{Computational Experiments:}

\medskip

\noindent
For the random realization of triangulated
orientable surfaces with $10$ vertices
we chose $k=2^{15}=32768$ as the side length 
of the cube.

\smallskip

\begin{itemize}
\item It took an average of about $700$ 
      tries to realize one of the $233$ triangulations of $S^2$ with
      $10$ vertices (in non-convex position).

      \vspace{1.5mm}

\item It took an average of about $418000$ 
      tries to realize one of the $2109$ triangulated tori with $10$
      vertices.

      \vspace{1.5mm}

\item For the $865$ vertex-minimal $10$-vertex triangulations of the 
      orientable surface of genus~$2$ we initially set a limit 
      of $200$ million tries for every example.
      Random realizations were found for about $\frac{1}{5}$
      of the $864$ triangulations. All other of the $864$ realizations
      were found by recycling of coordinates.
      The computations were run for $3$ months
      on ten Pentium~R $2.8$ GHz processors.

      \vspace{1.5mm}

\item It happened a few times that for a given triangulation
      realizations with identical coordinates
      were found by different processors after different numbers
      of tries. (Thus the length of cycles of the random generator plays a role.)
\end{itemize}

\smallskip

By successively rounding the coordinates of one of the randomly realized surfaces
it is in most cases relatively easy to obtain realizations of the respective surface 
with much smaller coordinates. 
(For triangulations of $S^2$ it is an open problem whether there are 
\emph{convex} realizations with small coordinates; cf.\  \cite[Ex.~4.16]{Ziegler1995}.)
An enumerative search for realizations with
small coordinates was carried out in \cite{HougardyLutzZelke2005apre} 
(this time, the computations were run for eight days on ten Pentium~R $2.8$ GHz processors)
and~\cite{HougardyLutzZelke2005bpre}; see these two references and
\cite{Lutz_PAGE} for explicit coordinates and visualizations.

\begin{thm} {\rm (Hougardy, Lutz, and Zelke~\cite{HougardyLutzZelke2005apre})}
All the $865$ vertex-minimal $10$-vertex triangulations of the orientable surface 
of genus $2$ have realizations in the \mbox{$(4\times 4\times 4)$}-cube, but 
cannot be realized (in general position) in the $(3\times 3\times 3)$-cube.
\end{thm}
A posteriori, the existence of small triangulations explains 
the successfulness of the described random search for realizations. 
\index{triangulated surfaces!realization!small coordinates}

With a more advanced simulated annealing approach,
realizations of surfaces can be obtained much easier:

\begin{thm} {\rm (Hougardy, Lutz, and Zelke~\cite{HougardyLutzZelke2006pre})}
All $20$ vertex-minimal $10$-vertex triangulations of the orientable surface 
of genus $3$ can be realized geometrically in~${\mathbb R}^3$.
\end{thm}

\subsection*{Acknowledgment} 

The author is grateful to Ewgenij Gawrilow
for a fast \texttt{polymake} \cite{polymake} implementation
of the random realization procedure.

\bibliography{}

\providecommand{\bysame}{\leavevmode\hbox to3em{\hrulefill}\thinspace}
\providecommand{\MR}{\relax\ifhmode\unskip\space\fi MR }
\providecommand{\MRhref}[2]{%
  \href{http://www.ams.org/mathscinet-getitem?mr=#1}{#2}
}
\providecommand{\href}[2]{#2}
\begin{thebibliography}{10}

\bibitem{Altshuler1974}
A.~Altshuler, \emph{Combinatorial $3$-manifolds with few vertices}, J.\ Comb.\
  Theory, Ser.\ A\/ \textbf{16} (1974), 165--173.

\bibitem{AltshulerBokowskiSchuchert1996}
A.~Altshuler, J.~Bokowski, and P.~Schuchert, \emph{Neighborly $2$-manifolds
  with $12$~vertices}, J.\ Comb.\ Theory, Ser.\ A\/ \textbf{75} (1996),
  148--162.

\bibitem{AltshulerSteinberg1973}
A.~Altshuler and L.~Steinberg, \emph{Neighborly $4$-polytopes with
  $9$~vertices}, J.\ Comb.\ Theory, Ser.\ A\/ \textbf{15} (1973), 270--287.

\bibitem{AltshulerSteinberg1976}
\bysame, \emph{An enumeration of combinatorial $3$-mani\-folds with nine
  vertices}, Discrete Math. \textbf{16} (1976), 91--108.

\bibitem{Barnette1982b}
D.~Barnette, \emph{Generating the triangulations of the projective plane}, J.\
  Comb.\ Theory, Ser.\ B\/ \textbf{33} (1982), 222--230.

\bibitem{BarnetteEdelson1988}
D.~W. Barnette and A.~L. Edelson, \emph{All $2$-manifolds have finitely many
  minimal triangulations}, Isr.\ J.\ Math. \textbf{67} (1988), 123--128.

\bibitem{Bokowski2006pre}
J.~Bokowski, \emph{On heuristic methods for finding realizations of surfaces},
  preprint, 2006, 6 pages.

\bibitem{BokowskiBrehm1987}
J.~Bokowski and U.~Brehm, \emph{A new polyhedron of genus $3$ with $10$
  vertices}, Intuitive Geometry, {\rm Internat.\ Conf.\ on Intuitive Geometry,
  Si\'ofok, Hungary, 1985} (K.~B\"or\"oczky and G.~{Fejes T\'oth},
  eds.), Colloquia Mathematica Societatis J\'anos Bolyai, vol.~48,
  North-Holland, Amsterdam, 1987, pp.~105--116.

\bibitem{BokowskiBrehm1989}
\bysame, \emph{A polyhedron of genus $4$ with minimal number of vertices and
  maximal symmetry}, Geom.\ Dedicata\/ \textbf{29} (1989), 53--64.

\bibitem{BokowskiEggert1991}
J.~Bokowski and A.~Eggert, \emph{Toutes les r\'{e}alisations du tore de
  {M}oebius avec sept sommets/{A}ll realizations of {M}oebius' torus with $7$
  vertices}, Topologie Struct. \textbf{17} (1991), 59--78.

\bibitem{BokowskiGuedes_de_Oliveira2000}
J.~Bokowski and A.~{Guedes de Oliveira}, \emph{On the generation of oriented
  matroids}, Discrete Comput.\ Geom. \textbf{24} (2000), 197--208.

\bibitem{BokowskiSturmfels1989}
J.~Bokowski and B.~Sturmfels, \emph{Computational {S}ynthetic {G}eometry},
  Lecture Notes in Mathematics, vol. 1355, Springer-Verlag, Berlin, 1989.

\bibitem{BowenFisk1967}
R.~Bowen and S.~Fisk, \emph{Generation of triangulations of the sphere}, Math.\
  Comput. \textbf{21} (1967), 250--252.

\bibitem{Brehm1981}
U.~Brehm, \emph{Polyeder mit zehn {E}cken vom {G}eschlecht drei}, Geom.\
  Dedicata\/ \textbf{11} (1981), 119--124.

\bibitem{Brehm1987b}
\bysame, \emph{A maximally symmetric polyhedron of genus $3$ with $10$
  vertices}, Mathematika\/ \textbf{34} (1987), 237--242.

\bibitem{plantri}
G.~Brinkmann and B.~McKay, \emph{\texttt{plantri}: a program for generating
  planar triangulations and planar cubic graphs},
  \url{http://cs.anu.edu.au/people/bdm/plantri/}, 1996--2001, Version 4.1.

\bibitem{Brueckner1897}
M.~Br\"uckner, \emph{Geschichtliche {B}emerkungen zur {A}ufz\"ahlung der
  {V}ielflache}, Pr.\ Realgymn.\ Zwickau. \textbf{578} (1897).

\bibitem{Brueckner1931}
\bysame, \emph{{\"Uber} die {A}nzahl $\psi (n)$ der allgemeinen {V}ielflache},
  Atti Congresso Bologna\/ \textbf{4} (1931), 5--11.

\bibitem{Csaszar1949}
A.~{Cs\'asz\'ar}, \emph{A polyhedron without diagonals}, Acta Sci.\ Math.,
  Szeged\/ \textbf{13} (1949--1950), 140--142.

\bibitem{Datta1999}
B.~Datta, \emph{Two dimensional weak pseudomanifolds on seven vertices}, Bol.\
  Soc.\ Mat.\ Mex., III.\ Ser. \textbf{5} (1999), 419--426.

\bibitem{DattaNilakantan2002}
B.~Datta and N.~Nilakantan, \emph{Two-dimensional weak pseudomanifolds on eight
  vertices}, Proc.\ Indian Acad.\ Sci.\ (Math.\ Sci.) \textbf{112} (2002),
  257--281.

\bibitem{Duke1970}
R.~A. Duke, \emph{Geometric embedding of complexes}, Am.\ Math.\ Mon.
  \textbf{77} (1970), 597--603.

\bibitem{Fendrich2003}
S.~Fendrich, \emph{Methoden zur {E}rzeugung und {R}ealisierung von
  triangulierten kombinatorischen $2$-{M}annigfaltigkeiten}, Diplomarbeit,
  Technische Universit\"at Darmstadt, 2003, 56 pages.

\bibitem{Gardner1975b}
M.~Gardner, \emph{Mathematical {G}ames. {O}n the remarkable {C}s\'asz\'ar
  polyhedron and its applications in problem solving}, Scientific American
  \textbf{232} (1975), no.~5, 102--107.

\bibitem{polymake}
E.~Gawrilow and M.~Joswig, \emph{\texttt{polymake}, {V}ersion 2.2, 1997--2006,
  with contributions by {T}.~{S}chr\"oder and {N}.~{W}itte},
  \url{http://www.math.tu-berlin.de/polymake}.

\bibitem{Grace1965}
D.~W. Grace, \emph{Computer search for non-isomorphic convex polyhedra}, Report
  CS 15, Computer Science Department, Stanford University, 1965.

\bibitem{Gruenbaum1967}
B.~Gr\"unbaum, \emph{Convex {P}olytopes}, Pure and Applied Mathematics,
  vol.~16, Interscience Publishers, London, 1967, second edition (V.~Kaibel,
  V.~Klee, and G.~M.~Ziegler, eds.), Graduate Texts in Mathematics,
  vol.~221, Springer-Verlag, New York, NY, 2003.

\bibitem{Heawood1890}
P.~J. Heawood, \emph{Map-colour theorem}, Quart.\ J.\ Pure Appl.\ Math.
  \textbf{24} (1890), 332--338.

\bibitem{HougardyLutzZelke2005apre}
S.~Hougardy, F.~H. Lutz, and M.~Zelke\mbox{}\mbox{}, \emph{Polyhedra of genus
  $2$ with $10$ vertices and minimal coordinates}, \url{arXiv:math.MG/0507592},
  2005, 3 pages.

\bibitem{HougardyLutzZelke2005bpre}
\bysame, \emph{Polyhedral tori with minimal coordinates}, in preparation.

\bibitem{HougardyLutzZelke2006pre}
\bysame, \emph{Surface realization with the intersection edge functional}, in
  preparation.

\bibitem{JungermanRingel1980}
M.~Jungerman and G.~Ringel, \emph{Minimal triangulations on orientable
  surfaces}, Acta Math. \textbf{145} (1980), 121--154.

\bibitem{KoehlerLutz2005pre}
E.~G. K\"ohler and F.~H. Lutz\mbox{}\mbox{}, \emph{Triangulated {M}anifolds
  with {F}ew {V}ertices: {V}ertex-{T}ransitive {T}riangulations {I}},
  \url{arXiv:math.GT/0506520}, 2005, 74 pages.

\bibitem{KuehnelLassmann1985-di}
W.~K\"uhnel and G.~Lassmann, \emph{Neighborly combinatorial $3$-manifolds with
  dihedral automorphism group}, Isr.\ J.\ Math. \textbf{52} (1985), 147--166.

\bibitem{Lavrenchenko1990}
S.~A. Lavrenchenko, \emph{Irreducible triangulations of the torus}, J.\ Sov.\
  Math. \textbf{51} (1990), 2537--2543, translation from Ukr.\ Geom.\ Sb.\
  \textbf{30} (1987), 52--62.

\bibitem{LawrencenkoNegami1997}
S.~[A.]~Lawrencenko and S.~Negami, \emph{Irreducible triangulations of the {K}lein
  bottle}, J.\ Comb.\ Theory, Ser.\ B\/ \textbf{70} (1997), 265--291.

\bibitem{Lutz2002b}
F.~H. Lutz, \emph{Cs\'asz\'ar's torus}, Electronic Geometry Models \textbf{{\rm
  No.\ 2001.02.069}} (2002), \url{http://www.eg-models.de/2001.02.069}.

\bibitem{Lutz_PAGE}
F.~H. Lutz\mbox{}, \emph{The {M}anifold {P}age, 1999--2006},
  \url{http://www.math.tu-berlin.de/diskregeom/stellar/}.

\bibitem{nauty}
B.~D. McKay, \emph{\texttt{nauty}, {V}ersion 2.2},
  \url{http://cs.anu.edu.au/people/bdm/nauty/}, 1994--2003.

\bibitem{Moebius1886}
A.~F. M\"obius, \emph{Mittheilungen aus {M}\"obius' {N}achlass: {I}.\ {Z}ur
  {T}heorie der {P}oly\"eder und der {E}lementarverwandtschaft}, Gesammelte
  Werke II (F.~Klein, ed.), Verlag von S.~Hirzel, Leipzig, 1886, pp.~515--559.

\bibitem{Ringel1955}
G.~Ringel, \emph{Wie man die geschlossenen nichtorientierbaren {F}l\"achen in
  m\"oglichst wenig {D}reiecke zerlegen kann}, Math.\ Ann. \textbf{130} (1955),
  317--326.

\bibitem{Royle_url}
G.~F. Royle, \emph{Number of planar triangulations},
  \url{http://www.csse.uwa.edu.au/~gordon/remote/planar/index.html#pts}.

\bibitem{Schewe2006pre}
L.~Schewe, 2006, work in progress.

\bibitem{Steinitz1922}
E.~Steinitz, \emph{Polyeder und {R}aumeinteilungen}, Encyklop\"adie der
  mathematischen Wissenschaften mit Einschluss ihrer Anwendungen, Dritter Band:
  Geometrie, {\rm III.1.2., Heft~9} (W.~Fr. Meyer and H.~Mohrmann, eds.),
  B.~G.~Teubner, Leipzig, 1922, pp.~1--139.

\bibitem{SteinitzRademacher1934}
E.~Steinitz and H.~Rademacher, \emph{Vorlesungen \"uber die {T}heorie der
  {P}olyeder unter {E}inschlu{\ss} der {E}lemente der {T}opologie}, Grundlehren
  der mathematischen Wissenschaften, vol.~41, Springer-Verlag, Berlin, 1934,
  reprint, 1976.

\bibitem{Sulanke2004pre}
T.~Sulanke, \emph{Note on the irreducible triangulations of the {K}lein bottle},
  \url{arXiv:math.CO/0407008}, 2004, 9~pages.

\bibitem{Sulanke2005apre}
\bysame, \emph{Generating triangulations of surfaces\mbox{}}, preprint,
  2005, 10 pages.

\bibitem{Sulanke2005bpre}
\bysame, \emph{Irreducible triangulations of {$S_2$}, {$N_3$}, and {$N_4$}},
  preprint, 2005, 8~pages.

\bibitem{Sulanke2005cpre}
\bysame, \emph{Source for \texttt{surftri} and lists of irreducible
  triangulations},
  \url{http://hep.physics.indiana.edu/~tsulanke/graphs/surftri/}, 2005,
  {V}ersion 0.96.

\bibitem{Tutte1962}
W.~T. Tutte, \emph{A census of planar triangulations}, Can.\ J.\ Math.
  \textbf{14} (1962), 21--38.

\bibitem{Ziegler1995}
G.~M. Ziegler, \emph{Lectures on {P}olytopes}, Graduate Texts in Mathematics,
  vol. 152, Springer-Verlag, New York, NY, 1995, revised edition, 1998.

\end{thebibliography}

\end{document}